\documentstyle[12pt]{article}
\newcommand{\nc}{\newcommand}
\nc{\al}{\alpha}
\nc{\ald}{{\dot \al}}
\nc{\ba}{\beta_\al}
\nc{\bb}{\beta_\beta}
\nc{\ga}{\g^\al}
\nc{\gb}{\g^\beta}
\nc{\db}{\pa_\beta}
\nc{\dtb}{\delta_\theta^\beta}
\nc{\dab}{{\delta_\al}^\beta}
\nc{\vmab}{V_{-\al}^\beta}
\nc{\vab}{V_\al^\beta}
\nc{\vib}{V_i^\beta}
\nc{\g}{\gamma}
\nc{\G}{\Gamma}
\nc{\D}{\Delta}
\nc{\paj}{P_{-\al}^j}
\nc{\la}{\lambda}
\nc{\La}{\Lambda}
\nc{\var}{\varphi}
\nc{\kvt}{\sqrt{t}}
\nc{\hn}{h^\vee}
\nc{\kn}{k^\vee}
\nc{\pa}{\partial}
\nc{\nn}{\nonumber \\ }
\nc{\hf}{\frac{1}{2}}         
\nc{\dz}{\frac{dz}{2\pi i}}
\nc{\fabc}{{f_{ab}}^c}
\nc{\binomial}[2]{\left (\begin{array}{c} {#1}\\ {#2} \end{array}
\right )}
\nc{\cpp}{{C_+}^+}
\nc{\cnp}{{C_0}^+}
\nc{\cmp}{{C_-}^+}
\nc{\cmn}{{C_-}^0}
\nc{\cmm}{{C_-}^-}
\nc{\vpp}{{V_+}^+}
\nc{\ben}{\begin{equation}}
\nc{\een}{\end{equation}}
\nc{\bea}{\begin{eqnarray}}
\nc{\eea}{\end{eqnarray}}
\nc{\bra}[1]{\langle {#1}|}
\nc{\ket}[1]{|{#1}\rangle}

\nc{\C}{\mbox{\hspace{1.24mm}\rule{0.2mm}{2.5mm}\hspace{-2.7mm} C}}
\nc{\Nat}{\mbox{\hspace{.04mm}\rule{0.2mm}{2.8mm}\hspace{-1.5mm} N}}

\nc{\spa}{\hspace{1 cm},\hspace{1 cm}}
\nc{\vs}{\vspace}
\nc{\NP}[1]{Nucl.\ Phys.\ {\bf #1}}
\nc{\PL}[1]{Phys.\ Lett.\ {\bf #1}}
\nc{\CMP}[1]{Commun.\ Math.\ Phys.\ {\bf #1}}
\nc{\PR}[1]{Phys.\ Rev.\ {\bf #1}}
\nc{\PRL}[1]{Phys.\ Rev.\ Lett.\ {\bf #1}}
\nc{\PTP}[1]{Prog.\ Theor.\ Phys.\ {\bf #1}}
\nc{\PTPS}[1]{Prog.\ Theor.\ Phys.\ Suppl.\ {\bf #1}}
\nc{\MPL}[1]{Mod.\ Phys.\ Lett.\ {\bf #1}}
\nc{\IJMP}[1]{Int.\ Jour.\ Mod.\ Phys.\ {\bf #1}}
\nc{\IM}[1]{Invent.\ Math.\ {\bf #1}}
\nc{\SJNP}[1]{Sov. J. Nucl. Phys.\ {\bf #1}}

\begin{document}

\topmargin -5mm
\oddsidemargin 5mm

\begin{titlepage}
\setcounter{page}{0}
\begin{flushright}
December 1997
\end{flushright}

\vs{8mm}
\begin{center}
{\Large Screening Current Representation of}\\[.2cm]
{\Large Quantum Superalgebras}

\vs{8mm}
{\large J{\o}rgen Rasmussen}\footnote{e-mail address: 
jorgen@celfi.phys.univ-tours.fr}\\[.2cm]
{\em Laboratoire de Math\'{e}mathiques et Physique Th\'{e}orique,}\\
{\em Universit\'{e} de Tours, F-37200 Tours, France}

\end{center}

\vs{8mm}
\centerline{{\bf{Abstract}}}
\noindent
In this letter a screening current or contour representation
is given of certain quantum superalgebras. The Gomez-Sierra construction
of quantum groups in conformal field theory is generalized to cover
superalgebras and illustrated using recent results on screening
currents in affine current superalgebra.\\[.4cm]

\end{titlepage}
\newpage
\renewcommand{\thefootnote}{\arabic{footnote}}
\setcounter{footnote}{0}

\section{Introduction}

Ever since their appearances \cite{FKRS,DJ} quantum groups have played an 
increasingly significant role in both mathematics and theoretical
physics. An important reason is due to their property of
bridging various areas of mathematical physics. Also quantum
supergroups have attracted a great deal of attention, see e.g. \cite{KR}
and the review \cite{CK}.

In the work \cite{GS1} by Gomez and Sierra a particular
construction of $q$-deformed universal enveloping algebras was
introduced in the context of minimal models in conformal field theory.
It is based on vertex operators screened by integrated
screening currents, and by use of contour manipulations
the quantum group structure is revealed. Their result allowed for
an immediate generalization to affine current algebras \cite{RRR},
and in \cite{GS2} a general scheme for obtaining the Gomez-Sierra
construction is discussed. However, the generality is limited to 
purely bosonic cases. The object of the present letter is to generalize
their scheme to include the case of superalgebras. The particular example
of affine current superalgebras \cite{Ras2} is then discussed.

The remaining part of this letter is organized as follows. In Section 2
we discuss in a general setting the Gomez-Sierra construction in the
case of superalgebras, closely following the original works on the
purely bosonic cases \cite{GS1,RRR,GS2}. In Section 3 we then illustrate
the generalization by considering affine current superalgebras 
using recent results on those obtained by the present author in \cite{Ras2}.
Section 4 contains some concluding remarks, in particular some speculations
on a possible generalization by including also screening currents of
the second kind \cite{BO,Ras1,PRY4,Ras3}. The construction obtained 
in affine current algebra
\cite{RRR} and the one in this letter in affine current superalgebra
are based only on screening currents of the first kind. 

\section{Gomez-Sierra Construction for Superalgebras}

The discussion here is in the framework of conformal field theory. It is 
assumed that the conformal field theory under consideration admits
a Coulomb gas or Feigin-Fuks representation. 
This ensures the existence of a set of screening currents where a 
screening current has the property
that the singular part of the operator product expansion 
with a chiral algebra generator is a total 
derivative. In particular, the screening currents are conformally primary 
of weight 1. These properties ensure that integrated
screening currents (screening charges) may be inserted into correlators
in the conformal field theory
without altering the chiral (including the conformal) Ward identities. 
This in turn makes them very useful in construction of correlators, see e.g. 
\cite{DF,F,PRY1,Ras1}. 

Alternatively one may use the screening currents
in representing certain quantum algebras \cite{GS1,RRR,GS2}. The contours 
used in the latter case are different from and simpler than the ones
needed in the former case. 
In this section we shall generalize the latter construction of
purely bosonic quantum groups to quantum supergroups.

Let $\left\{S_j(w)\right\}_{j\in{\cal S}}$ be a set of screening currents 
and let $V_\la(z)$ denote a generic
vertex operator of weight $\la$. The parities of the fields are defined
by
\ben
 p(V_\la)=0\spa p(S_j)=\left\{\begin{array}{ll} 1&\ \ ,\ \ j\in\tau\\
                                                0&\ \ ,\ \  
                                                j\in {\cal S}\setminus
                                                \tau
                              \end{array}\right.
\een
where $\tau$ is a subset of ${\cal S}$. The construction of the quantum
superalgebra based on the following defining braiding properties
\bea
 V_\la(z)V_{\la'}(w)&=&e^{i\pi\Omega_{\la\la'}}V_{\la'}(w)V_\la(z)\nn
 S_j(z)V_\la(w)&=&e^{i\pi\Omega_{j\la}}V_\la(w)S_j(z)\nn
 S_i(z)S_j(w)&=&e^{i\pi\hat{\Omega}_{ij}}S_j(w)S_i(z)\nn
  &=&e^{i\pi(\Omega_{ij}+p(S_i)p(S_j))}S_j(w)S_i(z)
\label{braiding}
\eea
is equivalent to the one used in the bosonic case \cite{GS2}. Thus we 
will follow that but, however, use a slightly different notation. 
The matrices $\Omega_{ij}$ and $\Omega_{\la\la'}$ are assumed symmetric.

First one introduces for every $z$ a vector space ${\cal V}_{\la,z}$, 
spanned by screened vertex operators $U_{\la,I}(z)$
\ben
 U_{\la,I}(z)=\int_{C_1}dw_1S_{i_1}(w_1)...\int_{C_{|I|}}dw_{|I|}
  S_{i_{|I|}}(w_{|I|})V_\la(z)
\een
where $I=\left\{i_1,...,i_{|I|}\right\}$ is an ordered set of
labels of screening currents.
The contours all start at infinity, encircle counter-clockwise the 
branch point $z$ once, and return to infinity. They do not intersect, nor do
they cross the branch cut from $z$ to infinity. The ordering of the contours
is chosen such that $C_i$ lies inside $C_j$ for $i>j$. Such contours will be
denoted standard contours around $z$. The vertex operator itself belongs
to ${\cal V}_{\la,z}$, namely $V_\la(z)=U_{\la,\emptyset}(z)$. It is called
the highest weight vector, distinguished by the absence of contour
structure.

For every screening current $S_j$ one defines a quartet of operators 
$F_j,K_j^{\pm1}$ and $E_j$. 
The former is defined as the contour creating operator
\ben
 F_j\left(U_{\la,I}(z)\right)=\int_CdwS_j(w)U_{\la,I}(z)
\een
where $C$ is a standard contour around $z$ and encloses the ones in
$U_{\lambda,I}(z)$. 
By this definition $F_j$ inherits the parity of $S_j$.
Now, choosing the branch cuts inherent in ${\cal V}_{
\la_1,z_1}$ and ${\cal V}_{\la_2,z_2}$ non-intersecting, we may
define the co-multiplication $\D(F_j)$ on ${\cal V}_{
\la_1,z_1}\otimes{\cal V}_{\la_2,z_2}$ by
\bea
 \D(F_j)\left(U_{\la_1,I_1}(z_1)\otimes U_{\la_2,I_2}(z_2)\right)
  &=&\int_CdwS_j(w)U_{\la_1,I_1}(z_1)U_{\la_2,I_2}(z_2)\nn
 &=&F_j\left(U_{\la_1,I_1}(z_1)\right)U_{\la_2,I_2}(z_2)\nn
 &+&e^{i\pi\left(\sum_{i\in I_1}\hat{\Omega}_{ji}
  +\Omega_{j\la_1}\right)}     
  U_{\la_1,I_1}(z_1)F_j\left(U_{\la_2,I_2}(z_2)\right)
\label{quaF}
\eea
The contour $C$ is of the standard type wrt the {\em pair} $z_1$ and
$z_2$. By an obvious contour deformation it may be written as a sum
$C=C_1+C_2$ where $C_1$ ($C_2$) is a standard contour around $z_1$ ($z_2$).
Thereby one has obtained
\ben
 \D(F_j)=F_j\otimes1+K_j^{-1}\otimes F_j
\een
where $K_j$ (and $K_j^{-1}$) is defined by
\ben
 K_j\left(U_{\la,I}(z)\right)=e^{-i\pi\left(\sum_{i\in I}\Omega_{ji}
  +\Omega_{j\la}\right)}U_{\la,I}(z)
\een
Note that
\ben
 (A_1\otimes A_2)(U_1\otimes U_2)=(-1)^{p(A_2)p(U_1)}A_1(U_1)\otimes
  A_2(U_2)
\een
The co-multiplication of $K_j^{\pm1}$ is trivial
\ben
 \D(K_j^{\pm1})=K_j^{\pm1}\otimes K_j^{\pm1}
\een

The raising operator $E_j$ is related to the operator $\hat{E}_j$ (see
below (\ref{E})) defined implicitly by
\bea
 L_{-1}U_{\la,I}(z)&=&\left[ L_{-1},\prod_{i_l\in I}\int_{C_l}dw_l
  S_{i_l}(w_l)\right] V_\la(z)+\prod_{i_l\in I}\int_{C_l}dw_lS_{i_l}(w_l)
  L_{-1}V_\lambda(z) \nn
 &=&-\sum_{j\in {\cal S}}\left(q_j-q_j^{-1}\right)
  S_j(\infty)\hat{E}_j
  U_{\la,I}(z)+\prod_{i\in I}F_{i}L_{-1}V_\la(z) 
\eea
The last term is a screened descendant of $V_\la$. At this point 
the factor $(q_j-q_j^{-1})$ is merely a normalization constant. In the 
application in Section 3 on affine current superalgebras we shall
see that choosing $q_j\neq1$ appropriately (and related to 
$e^{i\pi\Omega_{jj}}$), we obtain a $q$-deformed Lie superalgebra. From 
\bea
 \left[ L_{-1},\prod_{i_l\in I}\int_{C_l}dw_lS_{i_l}(w_l)\right] V_\la(z)
  &=&-\sum_{i_l\in I}(q_{i_l}-q_{i_l}^{-1})S_{i_l}(\infty)
  \frac{1-e^{2\pi i\left(\sum_{i'>i_l}\hat{\Omega}_{i_li'}+\Omega_{i_l\la}
  \right)}}{q_{i_l}-q_{i_l}^{-1}}\nn
 &\cdot&e^{i\pi\sum_{i'<i_l}\hat{\Omega}_{i_li'}}
  \prod_{i_{l'}\in I\setminus\left\{i_l\right\}}\int_{C_{l'}}dw_{l'}
  S_{i_{l'}}(w_{l'})V_\la(z)
\eea
we obtain
\ben
 \hat{E}_jU_{\lambda,I}=\sum_{i_l\in I,i_l\sim j}\frac{1-e^{2\pi i\left(
  \sum_{i'>i_l}\Omega_{ji'}+\Omega_{j\la}\right)}}{q_j-q_j^{-1}}
  e^{i\pi\sum_{i'<i_l}\hat{\Omega}_{ji'}}U_{\la,I\setminus\left\{i_l\right\}}
\label{Ehatqua}
\een
Here we have used that $U_{\la,I\setminus\left\{i\right\}}=0$ for $i\not\in
I$. The restriction in the summation means $i_l\sim j$ if $S_{i_l}=S_j$.
(\ref{Ehatqua}) is a simple generalization of the analogous result for purely
bosonic generators \cite{Ras1}. Note that just like $F_j$, $\hat{E}_j$
inherits the parity of $S_j$.
In the special case of only one type of screening current, (\ref{Ehatqua})
reduces to ($U_{\la;|I|}=U_{\lambda,I}$)
\ben
 \hat{E}_jU_{\lambda;|I|}=\frac{1-e^{2\pi i\Omega_{j\la}}
  e^{i\pi(|I|-1)\hat{\Omega}_{jj}}}{q_j-q_j^{-1}}
  \left[|I|\right]_{e^{i\pi\hat{\Omega}_{jj}}}U_{\la;|I|-1}
\een
This is in agreement with \cite{GS2} in the purely bosonic case. 
The $q$-number is introduced as
\ben
 [a]_q=\frac{1-q^a}{1-q}
\een

The co-multiplication of $\hat{E}_j$ is defined implicitly using 
the trivial co-multiplication of $L_{-1}$ (which is
$\D(L_{-1})=L_{-1}\otimes1+1\otimes L_{-1}$)
\bea
 &&\D(L_{-1})\left(U_{\la_1,I_1}(z_1)\otimes U_{\la_2,I_2}(z_2)\right)\nn
 &=&\prod_{i\in I_1}F_i L_{-1}V_{\la_1}(z_1)\otimes U_{\la_2,I_2}(z_2)+
  U_{\la_1,I_1}(z_1)\otimes\prod_{i\in I_2}F_iL_{-1}V_{\la_2}(z_2)\nn
 &-&\sum_{j\in{\cal S}}(q_j-q_j^{-1})S_j(\infty)\D(\hat{E}_j)\left(
  U_{\la_1,I_1}(z_1)\otimes U_{\la_2,I_2}(z_2)\right)
\eea
and is found to be
\ben
 \D(\hat{E}_j)=\hat{E}_j\otimes1+K_j^{-1}\otimes\hat{E}_j
\een
After having redefined the raising operator as  
\ben
 E_j=K_j\hat{E}_j
\label{E}
\een
the co-multiplication reads
\ben
 \D(E_j)=E_j\otimes K_j+1\otimes E_j
\een

With the following definitions of an anti-pode $\g$ 
\ben
 \g(E_j)=-E_jK_j^{-1}\spa\g(K_j^{\pm1})=K_j^{\mp1}\spa \g(F_j)=-K_jF_j
\een
and a co-unit $\epsilon$ 
\ben
 \epsilon(E_j)=0\spa\epsilon(K_j^{\pm1})=1\spa \epsilon(F_j)=0
\een
it is straightforward to verify that the above is indeed a Hopf algebra.
A final task is to implement quasi-triangularity in the Hopf algebra,
which is done by defining the universal $R$ matrix to be the braiding matrix 
of two screened vertex operators. We refer to \cite{GS2} for details.

Let us finally summarize the results for the (anti-)commutation relations
\bea
 K_iK_j&=&K_jK_i\nn
 K_iE_j&=&e^{i\pi\Omega_{ij}}E_jK_i\nn
 K_iF_j&=&e^{-i\pi\Omega_{ij}}F_jK_i\nn
 E_iF_j-(-1)^{p(E_i)p(F_j)}F_jE_i&=&
  \delta_{ij}\frac{K_i-K_i^{-1}}{q_i-q_i^{-1}}
\label{commgen}
\eea

\section{$q$-deformed Enveloping Lie Superalgebras}

In this section we shall illustrate the above construction by considering
the recently obtained explicit results on screening currents in affine
current superalgebras \cite{Ras2}. Let us first review some of the
results in that work while fixing our notation. 

The rank of the underlying Lie superalgebra is $r$. Choose a set of simple 
roots $\left\{\al_i\right\}_{i=1,...,r}$ such that the subset $\tau$ of 
$\left\{1,...,r\right\}$ corresponds to the odd simple roots. The parities of
the Chevalley generators are then defined by
\bea
   p(H_i)&=&0\nn
   p(E_i)=p(F_i)&=&\left\{\begin{array}{ll}1\ \ ,&\ \ i\in\tau\\
                                            0\ \ ,&\ \ i\not\in\tau
                                            \end{array}
       \right.
\eea 
In general, the parity of a raising $E_\al$
or a lowering operator $F_\al$ is 0 (1) if the
corresponding root $\al$ is even (odd).   
For the purpose of definiteness we choose a standard convention where 
the Cartan matrix is given by
\ben
 A_{ij}=\left\{\begin{array}{ll} \frac{2\al_i\cdot\al_j}{\al_i^2}\ \ 
   &,\ \ \al_i^2
   \neq0\\
   \al_i\cdot\al_j\ \ &,\ \ \al_i^2=0\end{array}\right.
\label{Cartan}
\een

We shall determine the $\Omega$ matrices (\ref{braiding}) working in the 
framework of free field realizations \cite{Wak,Ras1,PRY4,Ras2}. 
The free field realization of an affine current superalgebra is based
on a bosonic scalar field $\var_i$ for every Cartan index $i$ and a pair of
free bosonic (fermionic) ghost fields $(\beta_\al,\g^\al)$ 
for every positive even (odd) root $\al$. Their defining operator product
expansions are
\bea
 \var_i(z)\var_j(w)&=&G_{ij}\ln(z-w)\nn
 \beta_\al(z)\g^{\al'}(w)&=&\frac{\delta^{\al'}_{\al}}{z-w}
\eea
$G_{ij}$ is the Cartan part of the Cartan-Killing form.
The affine current superalgebra is characterized by the level or central 
extension $k$ which may conveniently be renormalized using the dual
Coxeter number $h^\vee$ by defining the parameter $t$ to be
\ben
 t=k+h^\vee
\een
Since we don't need the explicit free field realization of the affine currents
we shall only review the result for the screening currents \cite{Ras2} which
may be written in terms of the matrix $C$ defined by
\ben
 C_a^b(\g(w))=-\g^\al(w){f_{\al,a}}^b
\een
We shall use the convention of implicit summation over positive roots.
${f_{a,b}}^c$ are the structure constants.
In \cite{Ras2} (we refer to that work for more details) it is found that the
screening currents are of the form
\ben
 S_j(w)=-:\left[B(-C(\g(w)))\right]_{\al_j}^{\al'}\beta_{\al'}(w):
 :e^{-\al_j(H\cdot\var(w))/\kvt}: 
\een
where 
\ben
 \pa\var_i(z)\al(H\cdot\var(w))=\frac{\al(H_i)}{z-w}
\een 
The function $B$ is the generating function for the Bernoulli numbers
\ben
 B(u)=\frac{u}{e^u-1}=\sum_{n\geq0}\frac{B_n}{n!}u^n
\een
while the vertex operators are simply given by
\ben
 V_\la(z)=:e^{\la\cdot\var(z)/\kvt}:
\een
It is now straightforward to verify that the $\Omega$ matrices are given by
\bea
 \Omega_{\la\la'}&=&\la\cdot\la'/t\nn
 \Omega_{i\la}&=&-\al_i\cdot\la/t\nn
 \Omega_{ij}&=&\al_i\cdot\al_j/t
\label{braidsup}
\eea

After introducing the deformation parameter
\ben
 q=e^{i\pi/t}
\een
and the symmetrizing parameters $D_i$ (which symmetrize
the Cartan matrix, $D_iA_{ij}$ is symmetric) 
\ben
 D_i=\left\{\begin{array}{ll}\al_i^2/2\ \ &,\ \ \al_i^2\neq0\\
   1&,\ \ \al_i^2=0\end{array}\right.
\label{symm}
\een
we choose the deformation parameters $q_i$ of the preceding section as 
\ben
 q_i=q^{D_i}
\label{q}
\een
Using these and the results (\ref{braidsup}) for the $\Omega$ matrices, it is
immediately found that the (anti-)commutator relations (\ref{commgen}) may
be written as
\bea
 K_iK_j&=&K_jK_i\nn
 K_iE_j&=&q_i^{A_{ij}}E_jK_i\nn
 K_iF_j&=&q_i^{-A_{ij}}F_jK_i\nn
 E_iF_j-(-1)^{p(E_i)p(F_j)}F_jE_i&=&
  \delta_{ij}\frac{K_i-K_i^{-1}}{q_i-q_i^{-1}}
\eea

Now the question arises as to what quantum superalgebra the above construction
corresponds. Owing to its obvious resemblance, the natural proposal would 
be that it is the $q$-deformed universal enveloping algebra of the underlying 
Lie superalgebra {\bf g}, namely\footnote{Admissible representations 
\cite{KK,KW} are characterized by the parameter $t$ being rational. This 
means that in those cases the construction corresponds to a $q$-deformation 
with $q$ a root of unity.} ${\cal U}_{e^{i\pi/t}}$({\bf g}). 
However, as it has been observed by Scheunert 
\cite{Sch}, the problem of presenting Lie superalgebras (and thus
their quantum deformations) by Serre type relations is not a trivial
one. It occurs that one needs supplementary relations not naively defined
by the Cartan matrix. We leave it for future investigations to determine
what set of Serre type relations the above construction respects.

\section{Conclusion}

In this letter we have generalized the Gomez-Sierra construction of quantum
groups to the case of supergroups and illustrated the generalization
using recent results on screening currents in affine current
superalgebra. 

Another possible generalization is to include screening
currents of the second kind. Based on steps already taken along those
lines \cite{Ras1} we anticipate that in the purely bosonic case
one should encounter a semi-direct sum of two $q$-deformed universal
enveloping algebras of the underlying Lie algebra and of the transposed
of that, respectively, with the two $q$'s being related as $q'=q^{t^2}$.  
However, as it has been demonstrated in \cite{Ras3}, 
screening currents of the second kind are very complicated objects.
We intend to come back elsewhere with more discussions on screening
currents of the second kind and of the associated quantum groups.
\\[.2cm]
{\bf Acknowledgement}\\[.2cm]
The author is financially supported by the Danish Natural Science Research
Council, contract no. 9700517.

\end{document}